\documentclass[a4paper,11pt]{amsart}
\usepackage[margin=3cm]{geometry}
\usepackage{amsfonts,amsmath,amssymb,amscd,graphicx}
\usepackage[colorlinks]{hyperref}
\usepackage{mathrsfs}

\usepackage{tensor}

\usepackage{graphicx}

\usepackage{caption, subcaption}

\usepackage[english]{babel}
\usepackage[utf8]{inputenc}
\usepackage{ucs}
\usepackage{float}%
\usepackage{enumitem}
\usepackage{graphicx}
\usepackage{tikz}
\usepackage{indentfirst}
\usepackage{fancyhdr}
\usepackage[all]{xy}
\usepackage{physics}
\usepackage{tikz-cd}
\usepackage{comment}

\theoremstyle{theorem}
\newtheorem{theorem}{Theorem}[section]
\newtheorem{corollary}[theorem]{Corollary}
\newtheorem{lemma}[theorem]{Lemma}

\newtheorem{proposition}[theorem]{Proposition}

\theoremstyle{definition}
\newtheorem{definition}[theorem]{Definition}
\newtheorem{remark}[theorem]{Remark}

\numberwithin{equation}{section}

\newenvironment{preuve}[1][]
{\vskip 2mm  \noindent\emph{\bf Proof#1. }}{$\Box$ \vskip 2mm}
 
\newcommand{\beq}{\begin{eqnarray}}
\newcommand{\eeq}{\end{eqnarray}}
\newcommand{\beqe}{\begin{eqnarray*}}
	\newcommand{\eeqe}{\end{eqnarray*}}

\DeclareMathOperator{\vol}{vol}
\DeclareMathOperator{\genus}{genus}

\DeclareMathOperator{\inj}{inj}

\newcommand{\Nn}{\mathbb{N}}
\newcommand{\R}{\mathbb{R}}

\newcommand{\C}{\mathbb{C}}

\let\epsilon=\varepsilon

\begin{document}

\title{Lower bound for the Cheeger constant of random complex curves}
%\date
\author{\sc Michele Ancona and Damien Gayet}
\maketitle
\begin{abstract}
In this paper, we provide a lower bound for the Cheeger constant and the spectral gap for random complex curves in $\C P^2$. The complex curve is endowed with the restriction of the ambient Fubini-Study metric, and the probability measure is the Gaussian measure induced by the $\mathscr{L}^2$-Hermitian product on the space of complex homogeneous polynomials
of degree $d$ in $3$ variables. The proof relies on our previous  bounds for the systole and the curvature of random complex curves, together with an isoperimetric inequality for small ovals on complex curves. More generally, we establish such lower bounds for random complex curves within complex projective manifolds.
\end{abstract}

\section{Introduction} 
This paper addresses the spectral properties of random complex curves in projective manifolds, specifically focusing on a lower bound for their spectral gap. We begin by presenting our main result for random plane curves, i.e., random complex curves in $\C P^2$. 

A degree $d$ complex curve is the zero locus of a degree $d$ homogeneous polynomial $P\in \C^{\mathrm{hom}}_d[X_0,X_1, X_2]$, so a degree $d$ random complex curve is defined by a random degree $d$ homogeneous polynomial. We equip $\C^{\mathrm{hom}}_d[X_0,X_1, X_2]$ with the only $U(3)$-invariant Hermitian product, that is, the scalar product that makes the family of monomials $$\left\{\sqrt{ \frac{(d+2)!}{2 i!j!k!}}X_0^{i}X_1^jX_2^{k}\right\}_{i+j+k =d}$$ an unitary basis. 
Remark that the previous Hermitian product associated  appears as the $\mathscr{L}^2$-product associated with the Fubini-Study metric.
The random polynomials we will consider in this paper are then of the form
\begin{equation}\label{mesure2} P= \sum_{i+j+k =d} a_{i,j,k}\sqrt{ \frac{(d+2)!}{2 i!j!k!}}
X_0^{i}X_1^jX_2^{k}, \end{equation}
 where $(a_I)_{|I|=d}$ are independent complex Gaussian variables $\mathcal{N}_\C(0,1)$. Equivalently, we equip $\C^{\mathrm{hom}}_d[X_0,X_1, X_2]$ with the Gaussian measure $\mu_d$ associated with the previous Hermitian product.
 
\begin{remark}\label{fubinistudy} A Hermitian product on  $\C^{\mathrm{hom}}_d[X_0,X_1, X_2]$ induces a Fubini-Study metric on the projectivization $\mathbb{P}(\C^{\mathrm{hom}}_d[X_0,X_1, X_2])$ and then a natural probability measure on it. Since any two proportional degree $d$ polynomials define the same plane curve we could have chosen to take a random point of $\mathbb{P}(\C^{\mathrm{hom}}_d[X_0,X_1, X_2])$ in order to define a random plane curve. As the Fubini-Study probability on $\mathbb{P}(\C^{\mathrm{hom}}_d[X_0,X_1, X_2])$ measure is nothing but the pushforward of Gaussian measure on $\C^{\mathrm{hom}}_d[X_0,X_1, X_2]$ with respect to the natural projection $\C^{\mathrm{hom}}_d[X_0,X_1, X_2]\rightarrow \mathbb{P}(\C^{\mathrm{hom}}_d[X_0,X_1, X_2])$, the two probabilistic models are equivalent.
\end{remark}
 
 The vanishing locus $Z(P)\subset \C P^2$ of $P$ is almost surely a smooth connected Riemann surface of genus $\frac{(d-1)(d-2)}{2}$, 
 which is naturally equipped with a Riemannian metric: the restriction $g_{\mathrm{FS}|Z(P)}$ of the Fubini-Study metric $g_{\mathrm{FS}}$ of $\C P^2$.
 \begin{remark} As, almost surely, the topology of $Z(P)$ only depends on the degree $d$,  we can view $Z(P)$ as a \emph{fixed} genus $\frac{(d-1)(d-2)}{2}$ surface equipped with a \emph{random} metric. 
 \end{remark}
 
 Recall the definition of the Cheeger constant $h(S,g)$ of a Riemannian surface $(S,g)$:
  \[
 h(S,g)=\inf \frac{\ell(\gamma)}{\min\{\vol(A),\vol(B)\}}
 \] 
 where the infimum is taken among all smooth curves $\gamma $ which separates $S$ into two disjoint and non-empty part $A$ and $B$ and where $\ell(\gamma)$ is the length of $\gamma$.
 
 The goal of the paper is to give a probabilistic lower bound of the Cheeger constant $h(Z(P))$ of $(Z(P),g_{\mathrm{FS}|Z(P)})$. Remark that there is no non-trivial deterministic lower bound for $h(Z(P))$, that is, for any $d\geq 2$, one has 
 \[
 \inf_{P\in \C^{\mathrm{hom}}_d[X_0,X_1, X_2]} h(Z_P)=0.
  \]
  Indeed, let $P_1$ and $P_2$ be two polynomials of degree $d_1$ and $d_2$, respectively, with $d_1+d_2=d$. Suppose that $Z(P_1)$ and $Z(P_2)$ are smooth and intersect transversally (both are generic conditions). Then the degree $d$ polynomial $P=P_1P_2$ is such that $Z(P)$ is a reducible complex curve with $d_1d_2$ double points, denoted by $x_1,\dots,x_{d_1d_2}$. Note that, by construction, $Z(P)\setminus\{x_1,\dots,x_{d_1d_2}\}$ is not connected. Now, for any small $\epsilon>0$, we can take a small generic perturbation $Q_\epsilon$ of $P$ such that $Z(Q_\epsilon)$ contains $d_1d_2$ small loops (these are the vanishing cycles collapsing to $x_1,\dots,x_{d_1d_2}$ when $Q_\epsilon$ tends to $P$) whose length is smaller than $\epsilon$, and which separate $Z({Q_\epsilon})$ into two parts whose volumes are close to the volumes of $Z(P_1)$ and $Z(P_2)$.

The main result of this paper is a probabilistic lower bound for $h(Z(P))$.

\begin{theorem}\label{PlaneCheeger} One has
\[ 
\mu_d\left[h(Z(P))\geq\frac{1}{d^{5}}\right] \underset{d\to \infty}{\to}1
\]
where $\mu_d$ denotes the Gaussian measure on the space $\C^{\mathrm{hom}}_d[X_0,X_1, X_2]$.
\end{theorem}
The Cheeger constant $h(S,g)$ of a Riemannian surface $(S,g)$ is related to the spectral gap $\lambda_1(S,g)$ by the well-known \cite{cheeger} Cheeger  inequality $\lambda_1\geq  \frac{h^{2}}{4}$.  
In particular Theorem \ref{PlaneCheeger}  implies:
\begin{corollary}\label{lambda} One has
\[ 
\mu_d\left[\lambda_1(Z(P))\geq\frac{1}{d^{10}}\right] \underset{d\to \infty}{\to}1
\]
where $\mu_d$ denotes the Gaussian measure on the space $\C^{\mathrm{hom}}_d[X_0,X_1, X_2]$.
\end{corollary}
\begin{remark}
	\begin{enumerate}
\item The previous results actually hold in the much more general setting of random complex curves within complex projective manifolds, as introduced in Section \ref{sec general}.
\item 
Note that a (much worse) probabilistic exponential  lower bound \[ 
\mu_d\left[\lambda_1(Z(P))\geq\exp(-d^{6})\right] \underset{d\to \infty}{\to}1
\]
was obtained by the authors in \cite{AG} using methods different from those employed in this paper. It is also worth noting that, according to \cite{bourguignon}, one always has $\lambda_1(Z(P))\leq 6$, for any $P\in \C^{\textrm{hom}}_d[X_0,X_1,X_2]$. Can a better probabilistic upper bound be obtained?
\item 
For genus $g$ hyperbolic surfaces, that is, for compact surfaces with constant Gaussian curvature equal to $-1$, one has $\lambda_1\leq a_g$, where $a_g$ is a quantity depending only on the genus of the surface, and such that $a_g\rightarrow 1/4$ as $g\rightarrow \infty$, see \cite{huber}.
Recall that the volume of degree $d$ complex curves in $\C P^2$ is $d$, and that their average curvature is asymptotic to $-d$. Hence, in order to compare these complex curves to the hyperbolic surfaces, we should multiply the Fubini--Study metric by $d$. In this case the previous upper bound $6$ becomes $6/d$, to be compared to $1/4+o(1)$.  Note that there always exist high genus hyperbolic surfaces with $\lambda_1$ close to $1/4$, see \cite{hide}.
\end{enumerate}
\end{remark}
 
\subsection{General setting}\label{sec general}
Let $X$ be a complex projective manifold of dimension $n\geq 1$ together with  an ample holomorphic line bundle  $L\to X$. We equip $L$ with an Hermitian metric $h$  with positive curvature $\omega$, 
that is, locally 
\[ \omega = \frac{1}{2i\pi}\partial \bar{\partial}\log \|s\|^2_h>0,
\]
where $s$ is any local non vanishing holomorphic section of $L$. We denote by $g_\omega= \omega (\cdot, i\cdot)$ the associated K\"ahler metric on $X$. Let $(E,h_E)\to X$ be a holomorphic vector bundle of rank equal to $n-1$ equipped with a Hermitian metric $h_E$. For $d$ large enough, the vanishing locus $Z(s)$ of a generic section $s\in H^0(X,E\otimes L^d)$ is a smooth complex curve  of genus  satisfying (see for instance~\cite{ishihara1999generalization}):
	$$ 2\genus(Z(s))-2 = (K_X+c_1(E\otimes L^d))c_{n-1}(E\otimes L^d), $$
	where $K_X$ denotes the canonical bundle of $X$. 
	In particular, \begin{equation}\label{genre} 
	\genus(Z(s))\underset{d\to \infty}{\sim} \frac12\int_X \omega^{n} d^{n}.
	\end{equation}
	
The space $H^0(X,E\otimes L^{d})$ of holomorphic sections of $E\otimes L^{\otimes d}$ can be equipped with the $\mathscr{L}^2$-Hermitian product 
\begin{equation}\label{hermprod}
(s_1,s_2)\in (H^0(X,E\otimes L^d))^2\mapsto \langle s,t\rangle = \int_X \langle s_1(x),s_2(x)\rangle_{h_d}\frac{\omega^n}{n!},
\end{equation}
where 
$h_d:= h_E\otimes h^d.$

This Hermitian product induces a Gaussian measure $\mu_d$ over $H^0(X,E\otimes L^d)$, 
that is for any Borelian $U\subset H^0(X,E\otimes L^d),$
\begin{equation}\label{mesure}
 \mu_d (U)=\int_{s\in U} e^{-\frac12 \|s\|^2} \frac{\mathrm{d}s}{(2\pi)^{N_d}},
\end{equation}
where $N_d$ denotes the complex dimension of $H^0(X,E\otimes L^d)$ and $\mathrm{d} s$ denotes the Lebesgue measure associated to the Hermitian product~(\ref{hermprod}). 

If $(S_i)_{i\in \{1, \cdots, N_d\}}$
denotes an orthonormal basis of this space, then 
$$ s =\sum_{i=1}^{N_d} a_i S_i$$
follows the law $\mu_d$ if the random complexes $a_i$ are i.i.d standard complex Gaussians. 

As is Remark \ref{fubinistudy}, for any event depending only on the vanishing locus $Z(s)$ of a section $s\in H^0(X,E\otimes L^d)$, the probability measure $\mu_d$ can be replaced by  the Fubini--Study measure on the linear system $\mathbb{P} (H^0(X,E\otimes L^d))$.

\begin{theorem}\label{main theorem} Let $X$ be a $n$-dimensional projective manifold, $(L,h)$ be a positive Hermitian line bundle and $(E,h_E)$ be a rank $n-1$ Hermitian line bundle.
Let $(a_d)_d$ be a sequence of positive real numbers converging to $0$. Then, there exists $C>0$ such that
 $$\mu_d\left[h(Z(s))\geq \frac{a_d }{Cd^{\frac{5n-1}{2}}\sqrt{\log d}} \right]\geq 1- C(\sqrt{a_d}+\frac{1}d).$$
\end{theorem}
Using Cheeger inequality \cite{cheeger}, we obtain:
\begin{corollary} Let $(a_d)_d$ be a sequence of positive real numbers converging to $0$. Under the hypotheses of Theorem \ref{main theorem}, there exists $C>0$ such that
 $$\mu_d\left[\lambda_1(Z(s))\geq \frac{a^2_d }{Cd^{5n-1}\log d} \right]\geq 1- C(\sqrt{a_d}+\frac{1}d).$$
\end{corollary}

\subsection{Idea of proof} 

 Let $s$ be a holomorphic section of $E\otimes L^d$ and consider the Riemann surface $Z(s)$. Given $\gamma$  a (not necessarily connected) real curve in $Z(s)$ that separates $Z(s)$ in two parts $A$ and $B$, the goal is to bound from below $$h(\gamma):=\frac{\ell(\gamma)}{\min\{\vol(A),\vol(B)\}}.$$ 

We consider two cases. 
      \begin{itemize}
 \item If every connected component of $\gamma$ is a \emph{small} contractible loop, then we estimate $h(\gamma)$ using an isoperimetric inequality. Here, “small" means that every such loop is contained in a ball of radius  $\left(\sup_{x\in Z(s)}|K(x)|\right)^{-\frac{1}{2}}$, where $K$ is the Gaussian curvature of $Z(s)$. Such isoperimetric inequality is the novelty of the paper. The precise statement is given by Corollary \ref{hgamma}.
 
\item  If there is at least one connected component of  $\gamma$ which is not a small contractible loop, then such connected component is either non-contractible, or it is a \emph{big} contractible loop. In the first case the length of such component must be bigger than the systole of $Z(s)$. In the second case, the length is at least $\left(\sup_{x\in Z(s)}|K(x)|\right)^{-\frac{1}{2}}$.
On the other hand, $\min\{\vol(A),\vol(B)\}$ is always smaller than $\vol(Z(s))$, which is a deterministic quantity computed in  Lemma \ref{volume}. 

 \end{itemize}
Theorem \ref{main theorem} then follows from  bounds of the curvature  and of the systole of a random complex curve $Z(s)$ obtained by the authors in \cite{AG}.\\

\noindent
{\bf Acknowledgments.} The second author thanks G\'erard Besson for a valuable discussion on this subject. 
The research leading to these results has received funding from the French Agence nationale de la ANR-20-CE40-0017 (Adyct).

 \section{Length of ovals, area of disks and proof of Theorem \ref{main theorem}}
In this section, we will assume the hypotheses of 
Theorem~\ref{main theorem}.  Let $s\in H^0(X,E\otimes L^d)$  be a global transverse holomorphic section of $E\otimes L^d$, and $Z(s)$ be the zero locus of $s$. In particular, $Z(s)$ is a smooth complex curve. 
Let us denote by $V_d$ the volume of $Z(s)$ with respect to the restriction of the K\"ahler metric $g_\omega$ of $X$.
%%%%%%%%%%%%%%
\begin{lemma}\label{volume}Under the hypotheses above, 
 $$V_d=\sum_{i=0}^{n-1}\left(\int_X c_{n-1-i}(E)c_1(L)^{i+1}\right)d^i.$$
\end{lemma}
\begin{preuve}  Wirtinger inequality says that $V_d=\int_{Z(s)}\omega$. The latter is a cohomological quantity equal to $\langle [Z(s)],[\omega] \rangle$, where $[Z(s)]$ is the fundamental class of $Z(s)$ and $[\omega]$ the class of $\omega$. On one hand, we have $[\omega]=c_1(L)$. On the other hand, the Poincaré dual of $[Z(s)]$ equals $$c_{n-1}(E\otimes L^d)=\sum_{i=0}^{n-1}c_{n-1-i}(E)c_1(L^d)^{i}=\sum_{i=0}^{n-1}(c_{n-1-i}(E)c_1(L)^{i})d^i,$$
see \cite[Example 3.2.2]{fulton}.
This implies  $$V_d=\sum_{i=0}^{n-1}\left(\int_X c_{n-1-i}(E)c_1(L)^{i+1}\right)d^i.$$
Hence the result.
\end{preuve}

Under the hypotheses above, let $\gamma$ be a finite union of disjoint simple closed curves in $Z(s)$. We say that $\gamma$ is \emph{separating} is $Z(s)\setminus \gamma$ is not connected. 
Let 
\begin{equation}\label{separating}
\mathcal{C}(s)=\{\gamma\subset Z(s), \gamma \text{ is separating}\}.
\end{equation}
For any $\gamma\in \mathcal{C}(s)$, we can write  $Z(s)\setminus \gamma= A\cup B$, with $A$ and $B$ non-empty open disjoint subsets. We stress that $A$ and $B$ do not need to be connected (in particular, the partition of $Z(s)\setminus \gamma$ into two disjoint set is not unique). For such $\gamma\in\mathcal{C}(s)$, we define 
\begin{equation}\label{hg}
h(\gamma):=\frac{\ell(\gamma)}{\min\{\vol(A),\vol(B)\}}. 
\end{equation}
 where $\ell(\gamma)$ is the length of $\gamma$ and where the minimum in taken among any pairs of non-empty open disjoint subsets $A$ and $B$ with $Z(s)\setminus \gamma= A\cup B$. Note that 
\begin{equation}\label{inegalitefacile}
\forall \gamma\in\mathcal{C}(s), \ h(\gamma)>\frac{2\ell(\gamma)}{V_d},
\end{equation}
where $V_d$ is the volume of $Z(s)$.

\begin{definition}[Oval] An \emph{oval} in $Z(s)$ is a simple closed curve which bounds a disk. If $Z(s)$ is a smooth complex curve of positive genus, and $\tau$  is an oval in $Z(s)$, we denote by $D_\tau$ the unique disk bounded by $\tau$.
\end{definition}
Remark that by~(\ref{genre}), for $d$ large enough, the genus of $Z(s)$ is always positive.

\begin{definition} We say that a simple closed curve $\tau$ in $Z(s)$ is contained in a ball of radius $r$ is there exists $x\in Z(s)$ such that $\tau \subset B(x,r)$, where $B(x,r)$  denotes the geodesic ball (with respect to the metric $g_{\omega|Z(s)}$) of radius $r$ centered at $x$.
\end{definition}
Remark that if a simple closed curve $\tau$ in $Z(s)$ is contained in a ball of radius $r>0$ with $r<\inj(Z(s))$, then   $\tau$ is necessarily an oval.

For any $r>0$, let 
\begin{equation}\label{belo}
\mathcal{O}_r(s)=\{\tau\subset Z(s),\  \tau \text{ is an oval contained in a ball of radius } r\}.
\end{equation}
The following is the key technical result of the section.
\begin{proposition}\label{isoperimetric} Under the hypotheses of Theorem~\ref{main theorem}, let $(a_d)_d$ be a sequence of positive real numbers converging to $0$. There exists $C>0$ such that
 $$\mu_d\left[s\in H^0(X, E\otimes L^d), \inf_{\tau\in \mathcal{O}_{r_d}(s)}\frac{\ell(\tau)^2}{\vol(D_\tau)}\geq 2\pi \right]\geq 1- C(\sqrt{a_d}+\frac{1}d)$$
 where  $r_d=\frac{a_d }{Cd^{{\frac{3}{4}(n+1)}}\log d}$.
\end{proposition}
\begin{preuve} Let $d\in \Nn^*$ and $s\in H^0(X,E\otimes L^d)$ be a generic section such that $Z(s)$ has positive genus. 
By the Gauss equation in K\"ahler manifolds (see for example \cite[Proposition 9.2]{jost2008riemannian} or \cite[Theorem 2.2]{AG2}), the Gauss curvature of a complex curve in $X$ is bounded from above by the holomorphic sectional curvature of $X$. In particular, there exists a constant $K_+\in \R$ such that the curvature of $Z(s)$ is bounded from above by  $K_+$.
Now,  by the isoperimetric inequality on surfaces (see for example \cite{fiala} or \cite[Corollary 1 (iii)]{Topping}), if  $\tau$ is an oval of $Z(s)$, then $$\ell{(\tau})^2\geq \vol(D_\tau)(4\pi -K_+\vol(D_\tau)).$$ In order to prove the the proposition, it is then enough to prove that with large probability, 
	$$\tau\in \mathcal{O}_{r_d}(s)\Rightarrow K_+\vol(D_\tau)\leq  2\pi.$$ 
We will actually prove that 
there exists a sequence $(w_d)_d$ converging to $0$ as $d\to \infty$,  such that with large probability,
$$\tau\in \mathcal{O}_{r_d}(s)\Rightarrow
\vol(D_\tau)\leq w_d.$$ To do this, recall that by \cite[Theorem 1.6]{AG}, 
there exists $C$ (depending only on  $(a_d)_d, X, (L,h)$ and $(E,h_E)$) such that,
 $$\mu_d\left[s\in H^0(X, E\otimes L^d), \inf_{x\in Z(s)} K(x)\geq -k_d \right]\geq 1- C(\sqrt{a_d}+\frac{1}d),$$
with $k_d:=\frac{C}{a_d}d^{\frac{3}{2}(n+1)}\log^2d.$
Let $s$ be a section such that the curvature of $Z(s)$ is bounded from below by $-k_d$.
By the Bishop--Gromov inequality \cite[Corollary 4, page 256]{bishop}, the volume of a ball of radius $r_d$ in $Z(s)$ is less than or equal to the volume of a ball of the same radius $r_d$ inside the hyperbolic plane of constant curvature $-k_d$. The latter equals $$\vol(B_{h_{-k_d}}(x,r_d))=k^{-1}_d\vol(B_{h_{-1}}(x,\sqrt{k_d}r_d))=k^{-1}_d\vol(B_{h_{-1}}(x,C''\sqrt{a_d})=:w_d,$$ 
where $h_{-K}$ denotes the hyperbolic metric of constant curvature equal to $-K$ and $C''$ is a constant  depending only on $(a_d), X, (L,h)$ and $(E,h_E)$.  Note that $w_d$ tends to zero as $d\to \infty$  because $$\forall r>0, \ \vol(B_{h_{-1}}(x,r)= 2\pi(\mathrm{cosh}(r)-1).$$ Hence, the result follows.
\end{preuve}

For any $r>0$, define 
	$$\mathcal{C}_{r}(s)=\{\gamma \in \mathcal{C}(s),\  \forall \tau \text { connected component of } \gamma,\  \tau \in \mathcal O_r(s)\},$$
where $\mathcal{C}(s)$ has been defined by (\ref{separating}) and $\mathcal O_r(s)$ by (\ref{belo}).

\begin{definition}
	Under the hypotheses of Theorem~\ref{main theorem}, let $d\in \Nn^*$ and $s\in H^0(X,E\otimes L^d)$ be a generic section such that $Z(s)$ has positive genus.
		\begin{itemize}
	\item 
	We say that two disjoint ovals in $Z(s)$ form an \emph{injective pair} if one of the ovals is contained in the disk bounded by the other one. 
\item 	A \emph{nest} $N$ in $Z(s)$ is a union of disjoint ovals with the property that each pair of such ovals is an injective pair.
\item 
Let $\gamma\in\mathcal{C}_{r}(s)$. A \emph{full nest $N$ of $\gamma$} is a nest $N$ contained in $\gamma$ with the property that for each oval $\tau$ of $\gamma\setminus N$, the curve $\tau\cup  N$ is not a nest.  
\item The outermost oval  of a full nest $N$ is called a  \emph{maximal oval} and denoted by $\gamma_N$, the other ovals are called \emph{inner ovals}.
\end{itemize}
\end{definition}
We emphasize that an oval $\tau$ of $\gamma$ such that $D_\tau$ does not contain any other ovals of $\gamma$  is a full nest and that any $\gamma$ in  $\mathcal{C}_{r}(s)$ is the union of its full nests. 
\begin{proposition}\label{hofovals} Let $M, r>0$ be two positive real numbers. Under the hypotheses of Theorem~\ref{main theorem}, for any $d\in \Nn^*$ and $s\in H^0(X,E\otimes L^d)$ such that $Z(s)$ is smooth and has positive genus, we have 
$${\inf_{\tau\in \mathcal{O}_{r}(s)}\frac{\ell(\tau)^2}{\vol(D_\tau)}\geq M^2}\Rightarrow \inf_{ \gamma\in \mathcal{C}_{r}(s)} h(\gamma)\geq \frac{M\sqrt{2}}{\sqrt{V_d}},$$
where $V_d=\vol(Z(s))$ and $h(\gamma)$ is defined by~(\ref{hg}).
\end{proposition}
\begin{preuve}  
Let $\gamma\in  \mathcal{C}_{r}(s)$ and write $Z(s)\setminus \gamma= A\cup B$ with $A$ and $B$ non-empty open disjoint subsets. Denote by $N_1,\dots,N_k$ the full nests of  $\gamma$. For any $i\in\{1,\dots,k\}$, let $\gamma_{N_i}$ be the maximal oval of $N_i$ and $D_{N_i}$ be the disk bounded by $\gamma_{N_i}$. 
Let $$D:= \cup_i D_{N_i}\subset Z(s).$$ 
Since $D$ is a union of disjoint discs, $D^c$ is connected, and since $\gamma = \cup_i N_i$, $D^c$ does not intersect $\gamma$.
	Consequently,
one of the two sets $A$ or $B$ (say $A$) contains $D^c$ while   for any $i\in\{1,\dots,k\}$, $B$ contains a non-empty subset of $D_{N_i}$. For any $i\in\{1,\dots,k\}$, let $B_i=B\cap D_{N_i}$, so that $B=\bigcup_{i=1}^kB_i$.
By hypothesis,  $$\ell(\gamma)\geq\sum_{i=1}^k\ell(\gamma_{N_i})\geq M \sum_{i=1}^k\sqrt{\vol(D_{N_i})}.$$  This implies that 
$$
 \ell(\gamma)\geq 
 M\sum_{i=1}^{k}\vol(B_i)^{\frac{1}{2}}\geq M\left(\sum_{i=1}^{k}\vol(B_i)\right)^{\frac{1}{2}}=M\vol(B)^{\frac{1}{2}}.$$
If $\vol(B)\leq V_d/2$, then $$h(\gamma)\geq  \frac{M\vol(B)^{\frac{1}{2}}}{\min\{V_d-\vol(B),\vol(B)\}}=
M\frac{\vol(B)^{\frac{1}{2}}}{\vol(B)}\geq\frac{M \sqrt 2 }{\sqrt{V_d}} .$$ Similarly, if $\vol(B)\geq V_d/2$, we have $$M\frac{\vol(B)^{\frac{1}{2}}}{\min\{V_d-\vol(B),\vol(B)\}}=M\frac{\vol(B)^{\frac{1}{2}}}{V_d-\vol(B)}\geq \frac{M\left(V_d/2\right)^{\frac{1}{2}}}{V_d/2}=\frac{M\sqrt{2}}{\sqrt{V_d}}.$$
Hence the result.
\end{preuve}
Propositions \ref{isoperimetric} and \ref{hofovals} directly imply the following result.
\begin{corollary}\label{hgamma} Under the hypotheses of Proposition~\ref{hofovals}, let $(a_d)_d$ be a sequence of positive real numbers converging to $0$. 
Then, there exists $C>0$ such that
$$\mu_d\left[\inf_{ \gamma\in \mathcal{C}_{r_d}(s)} h(\gamma)\geq \sqrt{\frac{8\pi}{V_d}} \right]\geq 1- C(\sqrt{a_d}+\frac{1}d)$$
where $r_d=\frac{a_d }{Cd^{{\frac{3}{4}(n+1)}}\log d}$ and $V_d=\vol(Z(s))$.
\end{corollary}

We are now able to prove Theorem \ref{main theorem}.

\begin{preuve}[ of Theorem \ref{main theorem}] 
Let $(a_d)_d$ be a sequence of positive real number converging to $0$.
 For any $C>0$, define $A_d(C)$ to be the subset of sections  $s\in H^0(X,E\otimes L^d)$ such that 
 $$ \forall \gamma \subset Z(s), \gamma  \text{ is  non-contractible simple closed curve } \Rightarrow \ell(\gamma) \geq \frac{a_d d^{-\frac{3n+1}{2}}}{C\sqrt{\log d}}.$$ By \cite[Theorem 1.5]{AG},  we can find $C_1>0$ and $d_1\in\mathbb{N}$ such that  $$\forall d\geq d_1, \ \mu_d[A_d(C_1)]\geq 1-C_1(a_d+\frac{1}{d}).$$
 For any $C>0$, define $B_d(C)$ to be the subset of sections  $s\in H^0(X,E\otimes L^d)$ such that $$\inf_{ \gamma\in \mathcal{C}_{r_d(C)}(s)} h(\gamma)\geq \sqrt{\frac{8\pi}{V_d}},$$  where $r_d(C):=\frac{a_d }{Cd^{{\frac{3}{4}(n+1)}}\log d}$ and $V_d=\vol(Z(s))$. By Corollary \ref{hgamma},  we can find $C_2>0$ and $d_2\in\mathbb{N}$ such that  $$\forall d\geq d_2, \ \mu_d[B_d(C_2)]\geq 1-C_2(\sqrt{a_d}+\frac{1}d).$$ 
 Define $D_d:=A_d(C_1)\cap  B_d(C_2)$ and $d_3=\max\{d_1,d_2\}$. Then that there exists $C_3$ such that  
 $$ \forall d\geq d_3, \ \mu_d[D_d]\geq 1-C_3(\sqrt{a_d}+\frac{1}d).$$ 
 Take $s\in D_d$ and let $\gamma\in\mathcal{C}(s)$ (that is, $\gamma$ is seperating, recall~(\ref{separating})). We have two possible cases: either $\gamma\in\mathcal{C}_{r_d(C_2)}(s)$ or $\gamma\in\mathcal{C}(s)\setminus \mathcal{C}_{r_d(C_2)}(s)$. 
  In the first case, by definition of $B_d(C_2)$, one has $$h(\gamma)\geq \frac{\sqrt{8\pi}}{\sqrt{V_d}}\geq \frac{C_4}{d^{\frac{n-1}{2}}},$$ where we used that,  by Lemma  \ref{volume},  there exists $C_5>0$ depending only on $X$, $L$ and $E$ such that $$\forall d\geq 1, \ V_d\geq C_5 d^{n-1}.$$
  Let us now consider the case $\gamma\in\mathcal{C}(s)\setminus \mathcal{C}_{r_d(C_2)}(s)$. This happens if one of the following possibilities occurs:  either $\gamma$ has a non-contractible component, or $\gamma$ has an oval that is not contained in any ball of radius $r_d(C_2)$.
  If $\gamma$ has a non-contractible component, then by definition of $A_d$ and by Equation \eqref{inegalitefacile}, we have $$h(\gamma)\geq \frac{2a_d d^{-\frac{3n+1}{2}}}{C\sqrt{\log d}V_d}\geq \frac{C_6 a_d}{\sqrt{\log d}d^{\frac{5n-1}{2}}}.$$  If $\gamma$ has an oval that is not contained in any ball of radius $r_d(C_2)$, then the length of such oval is bigger than $r_d(C_2)$. In particular, by Equation \eqref{inegalitefacile}, we have $$h(\gamma)\geq \frac{r_d(C_2)}{V_d}\geq \frac{C_7a_d }{d^{\frac{7n-1}{4}}\log d}.$$ Hence, the result follows.
\end{preuve}
 
\bibliographystyle{amsplain}
\bibliography{cheeger.bib}
\vspace{12mm}

\noindent
Michele Ancona \\
Laboratoire J.A. Dieudonn\'e\\
UMR CNRS 7351\\
Universit\'e C\^ote d'Azur, Parc Valrose\\
06108 Nice, Cedex 2, France\\

\noindent
Damien Gayet\\
 Univ. Grenoble Alpes, Institut Fourier \\
F-38000 Grenoble, France \\
CNRS UMR 5208  \\
CNRS, IF, F-38000 Grenoble, France

\end{document}